\documentclass[11pt]{article}
\usepackage[margin=1in]{geometry}
\usepackage{enumerate}
\usepackage{url}
\usepackage{cite}

\usepackage{amsmath,amsthm,amssymb}

\title{\textbf{Approximation of maps between real algebraic varieties}}
\author{Juliusz Banecki and Wojciech Kucharz}
\date{}

\theoremstyle{plain}
\newtheorem{theorem}{Theorem}[section]
\newtheorem{lemma}[theorem]{Lemma}
\newtheorem{proposition}[theorem]{Proposition}
\newtheorem{corollary}[theorem]{Corollary}

\theoremstyle{definition}
\newtheorem{definition}[theorem]{Definition}
\newtheorem{example}[theorem]{Example}
\newtheorem{observation}[theorem]{Observation}
\newtheorem{notation}[theorem]{Notation, terminology and convention}

\newcounter{list}
\counterwithin*{list}{theorem}
\newenvironment{named}
  {\begin{list}
     {(\arabic{section}.\arabic{theorem}.\arabic{list})} 
     {\usecounter{list}   
      \setlength{\leftmargin}{1em}} 
  }
  {\end{list}}

\begin{document}

\maketitle

\paragraph{Abstract.} A nonsingular real algebraic variety $Y$ is said to have the approximation property if for every nonsingular real algebraic variety $X$ the following holds: if $f:X\rightarrow Y$ is a $\mathcal{C}^{\infty}$ map that is homotopic to a regular map, then $f$ can be approximated in the $\mathcal{C}^{\infty}$ topology by regular maps. In this paper, we characterize the varieties $Y$ with the approximation property. We also characterize the varieties $Y$ with the approximation property combined with a suitable interpolation condition. Some of our results have variants concerning the regular approximation of continuous maps defined on (possibly singular) real algebraic varieties.




\section{Introduction and main results}
\label{section:1}

In real algebraic geometry, the problem of approximating continuous or $\mathcal{C}^{\infty}$ maps by algebraic morphisms has been studied in numerous works since at least the 1980s (see \cite{1, 2, 5,6,7,8,9,10,11,12,13,14,15,16,17,18,26,27,28,31,32,33,35,37,40,43,45,47,49,50,51,56,58} and the references therein). The publications \cite{2,17} present the latest developments in this area. The line of research initiated in \cite{17} will be continued here with the introduction of new ideas and techniques.

Throughout this paper, we use the term \emph{real algebraic variety} to mean a ringed space with structure sheaf of $\mathbb{R}$-algebras of $\mathbb{R}$-valued functions, which is isomorphic to a Zariski locally closed subset of real projective $n$-space $\mathbb{P}^n(\mathbb{R})$, for some $n$, endowed with the Zariski topology and the sheaf of regular functions. Morphisms of real algebraic varieties are called \emph{regular maps}. This terminology is compatible with that used in the monographs \cite{5,49}, which contain a detailed treatment of real algebraic geometry. 

Recall that each real algebraic variety in the sense used here is actually \emph{affine}, that is, isomorphic to an algebraic set in $\mathbb{R}^n$ for some $n$ \cite[Proposition 3.2.10 and Theorem 3.4.4]{5} or \cite[Proposition 1.3.11]{49}. This fact allows us to describe regular functions and regular maps in a simple way. To be precise, let $X\subset \mathbb{R}^n$ and $Y\subset \mathbb{R}^p$ be algebraic sets, and let $U\subset X$ be a Zariski open set. A function $\varphi:U\rightarrow\mathbb{R}$ is regular if and only if there exist two polynomial functions $P,Q:\mathbb{R}^n\rightarrow\mathbb{R}$ such that 

$$Q(x)\neq 0 \textnormal{ and } \varphi(x)=\frac{P(x)}{Q(x)} \textnormal{ for all } x\in U$$
\cite[p.62]{5} or \cite[p.14]{49}. A map $f=(f_1,\ldots, f_p):U\rightarrow Y\subset\mathbb{R}^p$ is regular if and only if its components $f_1$, \ldots, $f_p$ are regular functions. 

Each real algebraic variety carries also the Euclidean topology determined by the usual metric on $\mathbb{R}$. Unless explicitly stated otherwise, all topological notions relating to real algebraic varieties refer to the Euclidean topology.

Given two $\mathcal{C}^{\infty}$ manifolds $M$ and $N$, we denote by $\mathcal{C}^k(M,N)$ the space of all $\mathcal{C}^k$ maps from $M$ to $N$, endowed with the $\mathcal{C}^k$ topology, where $k$ is a nonnegative integer or $k=\infty$ (see \cite[pp.34,36]{30} or \cite[p.311]{57} for the definition of this topology and note that in \cite{30} it is called the weak $\mathcal{C}^k$ topology; the $\mathcal{C}^0$ topology is just the compact-open topology). From now on, $\mathcal{C}^{\infty}$ manifolds and $\mathcal{C}^{\infty}$ maps will be referred to as \emph{smooth manifolds} and \emph{smooth maps}, respectively.

Let $X$, $Y$ be nonsingular real algebraic varieties, and $U$ an open subset of $X$. We say that a $\mathcal{C}^k$ map $f:U\rightarrow Y$ can be \emph{approximated in the $\mathcal{C}^k$ topology by regular maps from $X$ to $Y$} if for every open neighborhood $\mathcal{U}\subset \mathcal{C}^k(U,Y)$ of $f$ there exists a regular map $g:X\rightarrow Y$ whose restriction $g|_U$ belongs to $\mathcal{U}$.

According to Nachbin's version of the Stone-Weierstrass approximation theorem given for example in \cite[Theorem 1.2.1]{60}, every $\mathcal{C}^k$ function from $U$ to $\mathbb{R}$ can be approximated in the $\mathcal{C}^k$ topology by regular functions.

The Stone-Weierstrass theorem does not carry over directly to $Y$-valued maps. We encounter the following phenomena:
\begin{itemize}
    \item In \cite{26} Ghiloni introduced the notion of generic real algebraic variety and proved that the set of regular maps from $X$ to $Y$ is nowhere dense in $\mathcal{C}^0(X,Y)$ and in $\mathcal{C}^{\infty}(X,Y)$ if $Y$ is generic.
    \item If a continuous map from $X$ to $Y$ can be approximated in the $\mathcal{C}^0$ topology by regular maps, then it is homotopic to a regular map. This is because there is a strong deformation retraction of $X$ onto some compact subset \cite[Corollary 9.3.7]{5}.
    \item By \cite[Theorems 1.1 and 1.1']{15}, if $Y$ has a compact connected component of positive dimension, then there exists a smooth map from a compact nonsingular real algebraic variety $Z$ to $Y$ that is not homotopic to any regular map.
    \item It can happen that a continuous map from $X$ to $Y$ is homotopic to a regular map, but cannot be approximated by regular maps. For example, every regular map from $\mathbb{P}^1(\mathbb{R})$ to any nonsingular cubic curve $C$ in $\mathbb{P}^2(\mathbb{R})$ is constant, so among continuous maps from $\mathbb{P}^1(\mathbb{R})$ to $C$ only those that are constant can be approximated by regular maps.
\end{itemize}

This brings us to the issues addressed in the present paper.

\begin{definition}
\label{def:1.1}
    Let $Y$ be a nonsingular real algebraic variety.
    \begin{enumerate}
        \item[(1)] (Approximation property or AP).  The variety $Y$ is said to have the \emph{approximation property} (AP for short) if for every nonsingular real algebraic variety $X$ the following holds: if $f:X\rightarrow Y$ is a smooth map that is homotopic to a regular map from $X$ to $Y$, then $f$ can be approximated in the $\mathcal{C}^{\infty}$ topology by regular maps from $X$ to $Y$.
        \item[(2)] (Generalized approximation property or GAP).  The variety $Y$ is said to have the \emph{generalized approximation property} (GAP for short) if for every nonsingular real algebraic variety $X$ and every open subset $U$ of $X$ the following holds: if $f:U\rightarrow Y$ is a smooth map for which there exists a regular map $f_0:X\rightarrow Y$ with $f_0|_U$ homotopic to $f$, then $f$ can be approximated in the $\mathcal{C}^{\infty}$ topology by regular maps from $X$ to $Y$. 
    \end{enumerate}
\end{definition}

In the study of the approximation of smooth maps by regular maps, the focus has often been on the AP of the target variety (see \cite{1,5,6,8,9,16,17,18,31,35,37,40,43} and the references therein). Determining which varieties have the AP
 is a difficult problem. General results in this direction have been obtained only recently \cite{17}. For example, in \cite{17} it is proved that for every positive integer $n$, the unit $n$-sphere has the AP, while a weaker version of such a result was previously proved only for certain specific values of $n$, namely $n=1,2$ or $4$ in \cite{6} and $n=3$ or $7$ in \cite{1}.

 It turns out that AP and GAP are equivalent properties. The nonsingular real algebraic varieties having these properties are characterized in Theorem \ref{thm:1.3} below.

 We will work with vector bundles, which are always $\mathbb{R}$-vector bundles. Let $Y$ be a real algebraic variety. Given a vector bundle $p:E\rightarrow Y$, with total space $E$ and bundle projection $p$, we sometimes refer to $E$ as a vector bundle over $Y$. If $y$ is a point in $Y$, we let $E_y:=p^{-1}(y)$ denote the fiber over $y$ and write $0_y$ for the zero vector in $E_y$. We call $Z(E)=\{ 0_y\in E :\, y\in Y \}$ the zero section of $E$. The general theory of algebraic vector bundles over real algebraic varieties can be found in \cite[Section 12.1]{5}.

Assuming that $Y$ is a nonsingular real algebraic variety, we denote the tangent bundle of $Y$ by $TY$ and the tangent space to $Y$ at a point $y\in Y$ by $T_yY$.

For any topological space $S$, let $\Omega(S)$ denote the collection of all relatively compact open subsets of $S$. 

\begin{definition}
\label{def:1.2}
    Let $Y$ be a nonsingular real algebraic variety.
    \begin{enumerate}
        \item[(1)] (Quasi-spray).  A \emph{quasi-spray} $(E,p,\{s_V\},\{\tau_V\})$ for $Y$ consists of an algebraic vector bundle $p:E\rightarrow Y$ over $Y$ and two collections $\{s_V\}$, $\{\tau_V\}$ indexed by $V\in \Omega(Y)$, where $s_V:E\rightarrow Y$ is a regular map and $\tau_V:V\rightarrow E$ is a smooth map such that $s_V(\tau_V(y))=y$ and $p(\tau_V(y))=y$ for all $y\in V$.
        \item[(2)] (Dominating quasi-spray).  A quasi-spray $(E,p,\{s_V\},\{\tau_V\})$ for $Y$ is said to be \emph{dominating} if for every $V\in \Omega(Y)$ and every $y\in V$ the regular map
        $$E_y\rightarrow Y,\quad v\mapsto s_V(v+\tau_V(y))$$
        is a submersion at $0_y\in E_y$.
        \item[(3)] (Quasi-malleable variety).  The variety $Y$ is called \emph{quasi-malleable} if it admits a dominating quasi-spray.
        \item[(4)] (Spray).  A \emph{spray} $(E,p,s)$ for $Y$ consists of an algebraic vector bundle $p:E\rightarrow Y$ over $Y$ and a regular map $s:E\rightarrow Y$ such that $s(0_y)=y$ for all $y\in Y$.
        \item[(5)] (Dominating spray).  A spray $(E,p,s)$ for $Y$ is said to be \emph{dominating} if for every $y\in Y$ the regular map 
        $$E_y\rightarrow Y,\quad v\mapsto s(v)$$
        is a submersion at $0_y\in E_y$.
        \item[(6)] (Malleable variety).  The variety $Y$ is called \emph{malleable} if it admits a dominating spray.
    \end{enumerate}
\end{definition}

The notions of quasi-spray, dominating quasi-spray and quasi-malleable variety are new, while the concepts of spray, dominating spray and malleable variety in Definition \ref{def:1.2} are identical with those in \cite[Definition 2.1]{17}. Each spray $(E,p,s)$ for $Y$ determines a quasi-spray $(E,p,\{s_V\},\{\tau_V\})$ for $Y$, where $s_V=s$ and $\tau_V(y)=0_y$ for all $V\in \Omega(Y)$ and all $y\in V$. It follows that every malleable variety is quasi-malleable. On the other hand, the question whether every quasi-malleable variety is malleable seems to be a difficult problem to settle.

The simplest malleable variety is $Y=\mathbb{R}^n$ which admits a dominating spray $(E,p,s)$, where $p:E=Y\times \mathbb{R}^n\rightarrow Y$ is the product vector bundle and $s:E\rightarrow Y$ is defined by $s(y,v)=y+v$.

Two nontrivial classes of malleable varieties are identified in \cite{17, 41}; see Examples \ref{ex:1.7} and \ref{ex:1.8} below for a summary. New constructions of malleable varieties are given in Theorem \ref{thm:1.9} and Corollary \ref{cor:1.10}

The following theorem is one of our main results.

\begin{theorem}[Varieties with the AP]
\label{thm:1.3}
    For a given nonsingular real algebraic variety $Y$, the following conditions are equivalent:
    \begin{enumerate}
        \item[(a)] The variety $Y$ has the generalized approximation property.
        \item[(b)] The variety $Y$ has the approximation property.
        \item[(c)] Every smooth map $\sigma:TY\rightarrow Y$, with $\sigma(0_y)=y$ for all $y\in Y$, can be approximated in the $\mathcal{C}^1$ topology by regular maps from $TY$ to $Y$.
        \item[(d)] Every smooth map $\varphi:Y\times \mathbb{R}\rightarrow Y$, with $\varphi(y,0)=y$ for all $y\in Y$, can be approximated in the $\mathcal{C}^1$ topology by regular maps from $Y\times \mathbb{R}$ to $Y$.
        \item[(e)] The variety $Y$ is quasi-malleable.
    \end{enumerate}
\end{theorem}

Since all continuous maps defined on a contractible space are null homotopic, the following is an immediate corollary to Theorem \ref{thm:1.3}.

\begin{corollary}
\label{cor:1.4}
    Let $Y$ be a quasi-malleable nonsingular real algebraic variety. Let $X$ be a nonsingular real algebraic variety, and $U$ a contractible open subset of $X$. Then every smooth map from $U$ to $Y$ can be approximated in the $\mathcal{C}^{\infty}$ topology by regular maps from $X$ to $Y$.
\end{corollary}

We move on to examine approximation in combination with interpolation.

\begin{definition}
\label{def:1.5}
    Let $Y$ be a nonsingular real algebraic variety.
    \begin{enumerate}
        \item[(1)] (Approximation-interpolation property or AIP).  The variety $Y$ is said to have the \emph{approximation-interpolation property} (AIP for short) if for every nonsingular real algebraic variety $X$ and every nonsingular Zariski closed subvariety $A$ of $X$ the following holds: if $f:X\rightarrow Y$ is a smooth map for which there exists a regular map $f_0:X\rightarrow Y$ satisfying $f_0|_A=f|_A$ and homotopic to $f$ relative to $A$, then $f$ can be approximated in the $\mathcal{C}^{\infty}$ topology by regular maps from $X$ to $Y$ whose restrictions to $A$ are equal to $f_0|_A$.
        \item[(2)] (Generalized approximation-interpolation property or GAIP).  The variety $Y$ is said to have the \emph{generalized approximation-interpolation property} (GAIP for short) if for every nonsingular real algebraic variety $X$, every open subset $U$ of $X$ and every nonsingular Zariski closed subvariety $A$ of $X$ the following holds: if $f:U\rightarrow Y$ is a smooth map for which there exists a regular map $f_0:X\rightarrow Y$ with $f_0|_{U\cap A}=f|_{U\cap A}$ and with $f_0|_U$ homotopic to $f$ relative to $U\cap A$, then $f$ can be approximated in the $\mathcal{C}^{\infty}$ topology by regular maps from $X$ to $Y$ whose restrictions to $A$ are equal to $f_0|_A$.
    \end{enumerate}
\end{definition}

We have the following variant of Theorem \ref{thm:1.3}.

\begin{theorem}[Varieties with the AIP]
\label{thm:1.6}
    For a given nonsingular real algebraic variety $Y$, the following conditions are equivalent:
    \begin{enumerate}
        \item[(a)] The variety $Y$ has the generalized approximation-interpolation property.
        \item[(b)] The variety $Y$ has the approximation-interpolation property.
        \item[(c)] Every smooth map $\sigma:TY\rightarrow Y$, with $\sigma(0_y)=y$ for all $y\in Y$, can be approximated in the $\mathcal{C}^1$ topology by regular maps from $TY$ to $Y$ whose restrictions to $Z(TY)$ are equal to $\sigma|_{Z(TY)}$.
        \item[(d)] Every smooth map $\varphi:Y\times \mathbb{R}\rightarrow Y$, with $\varphi(y,0)=y$ for all $y\in Y$, can be approximated in the $\mathcal{C}^1$ topology by regular maps from $Y\times \mathbb{R}$ to $Y$ whose restrictions to $Y\times \{0\}$ are equal to $\varphi|_{Y\times \{0\}}$.
        \item[(e)] The variety $Y$ is malleable.
    \end{enumerate}
\end{theorem}

If $Y$ is a malleable variety, then according to Theorem \ref{thm:5.1} in Section \ref{section:5} it actually has a stronger property than the GAIP.

Theorems \ref{thm:1.3} and \ref{thm:1.6} result in the following diagram:
\begin{center}
    \begin{tabular}{c c c c c}
    GAP & $\iff$ & AP & $\iff$ & quasi-malleable \\
    $\Uparrow$ &  & $\Uparrow$ &  & $\Uparrow$ \\
    GAIP & $\iff$ & AIP & $\iff$ & malleable
\end{tabular}
\end{center}
It seems unlikely that the vertical implications can be reversed.

The relevance of Theorems \ref{thm:1.3} and \ref{thm:1.6} is reinforced by interesting examples of malleable varieties.

\begin{example}[Homogeneous spaces]
\label{ex:1.7}
    Let $G$ be a linear real algebraic group, that is, a Zariski closed subgroup of the general linear group $GL_n(\mathbb{R})$ for some $n$. A \emph{$G$-space} is a real algebraic variety $Y$ on which $G$ acts, where the action $G\times Y\rightarrow Y$, $(a,y)\mapsto a\cdot y$ is a regular map. We say that a $G$-space $Y$ is \emph{good} if $Y$ is nonsingular and for every point $y\in Y$ the map $G\rightarrow Y$, $a\mapsto a\cdot y$ is a submersion at the identity element of $G$. Obviously, if $Y$ is \emph{homogeneous} for $G$ (that is, $G$ acts transitively on $Y$), then $Y$ is a good space. The group $G$ acts transitively on itself by left multiplication and is therefore a homogeneous space. By \cite[Proposition 2.8]{17}, every good space is malleable.
\end{example}

In particular, the unit $n$-sphere 
$$\mathbb{S}^n:=\{ (x_0,\ldots, x_n)\in\mathbb{R}^{n+1}:\, x_0^2+\ldots+x_n^2=1 \}$$
and real projective $n$-space $\mathbb{P}^n(\mathbb{R})$ are malleable varieties because they are homogeneous spaces for the orthogonal group $O_{n+1}(\mathbb{R})\subset GL_{n+1}(\mathbb{R})$.

\begin{example}[Open subvarieties]
\label{ex:1.8}
    If $Y$ is a malleable real algebraic variety, then every Zariski open subset of $Y$ is also a malleable variety \cite[Proposition 2.11]{41}. In particular, every Zariski open subset of $\mathbb{R}^n$ is a malleable variety. It is obvious that the malleable variety $\mathbb{R}^n\setminus \mathbb{S}^{n-1}$, where $n\geq 1$, is not homogeneous.
\end{example}

Our new result is the following local malleability criterion.

\begin{theorem}[Localization theorem]
\label{thm:1.9}
    Let $Y$ be a nonsingular real algebraic variety, each point of which has a Zariski open neighborhood that is a malleable variety. Then $Y$ is a malleable variety.
\end{theorem}

It is worthwhile to compare malleable varieties with more familiar ones, such as rational or unirational varieties.

Recall that an $n$-dimensional irreducible real algebraic variety is said to be \emph{rational} if it contains a nonempty Zariski open subset biregularly isomorphic to a Zariski open subset of $\mathbb{R}^n$. We extend this definition slightly by declaring that a real algebraic variety is \emph{rational} if each of its irreducible components is rational. If a homogeneous space for some linear real algebraic group has at least one rational irreducible component, then it is a rational variety.

It is an open question whether every linear real algebraic group is a rational variety. However, there are homogeneous spaces for some linear real algebraic groups that are not rational varieties. For the clarification of this last point we are indebted to Oliver Benoist and Olivier Wittenberg. Thus, in view of Example \ref{ex:1.7}, a malleable variety is not necessarily rational. It is unknown whether every rational nonsingular real algebraic variety is malleable.

An $n$-dimensional real algebraic variety $Y$ is said to be \emph{uniformly rational} if every point in $Y$ has a Zariski open neighborhood that is biregularly isomorphic to a Zariski open subset of $\mathbb{R}^n$. Clearly, every uniformly rational variety is nonsingular. An intriguing open question posed by Gromov is whether every rational nonsingular variety is uniformly rational (see \cite[p.885]{29} and \cite{19} for the discussion in the context of complex algebraic varieties).

Uniformly rational varieties are important in our theory.

\begin{corollary}[Uniformly rational implies malleable]
\label{cor:1.10}
    Every uniformly rational real algebraic variety is malleable.
\end{corollary}
\begin{proof}
    According to Example \ref{ex:1.8}, every Zariski open subset of $\mathbb{R}^n$ is a malleable variety, so we complete the proof using Theorem \ref{thm:1.9}.
\end{proof}

\begin{example}[Uniformly rational varieties]
\label{ex:1.11}
    Let $Y$ be an irreducible nonsingular real algebraic variety.
    \begin{enumerate}
        \item[(1)] If $Y$ is rational and dim$Y=1$, then $Y$ is uniformly rational because it is biregularly isomorphic to a Zariski open subset of $\mathbb{P}^1(\mathbb{R})$.
        \item[(2)] If $Y$ is rational and dim$Y=2$, then as described in \cite[$\mathsection$2.2]{48}, $Y$ is uniformly rational by Comessatti's theorem \cite[p.257]{21}, whose modern proofs are given in \cite[Theorem 30]{34} and \cite[p.137, Proposition 6.4]{54}.
        \item[(3)] If $Y$ is uniformly rational and $Z$ is a nonsingular Zariski closed subvariety of $Y$, then the blow-up of $Y$ with center $Z$ is a uniformly rational variety. The proof given in \cite[p.885]{29} and \cite{19} in a complex setting also works for real algebraic varieties.
    \end{enumerate}
\end{example}

A real algebraic variety $Y$ is said to be \emph{unirational} if there exists a regular map $\varphi:U\rightarrow Y$ defined on a Zariski open subset $U$ of $\mathbb{R}^n$, for some $n$, such that the image $\varphi(U)$ is Zariski dense in $Y$. In that case, $Y$ is irreducible and one can choose a regular map $\varphi$ as above with $n=\textnormal{dim}Y$.

Every irreducible quasi-malleable real algebraic variety $Y$ is unirational. This is because if $(E,p,\{s_V\},\{\tau_V\})$ is a dominating quasi-spray for $Y$, then for every point $y\in V$ the regular map 
$$E_y\rightarrow Y,\quad v\mapsto s_V(v+\tau_V(y))$$
is a submersion at $0_y\in E_y$, so its image is Zariski dense in $Y$.

In summary, the following implications hold for all irreducible nonsingular real algebraic varieties:

\begin{center}
    \begin{tabular}{c c c}
        uniformly rational & $\Rightarrow$ & malleable \\
        $\Downarrow$ &  & $\Downarrow$ \\
        rational &  & quasi-malleable \\
         &  & $\Downarrow$ \\
          &  & unirational
    \end{tabular}
\end{center}

The converse of the horizontal implication is generally false, while there are no known counterexamples to the converse of any of the vertical implications.

It is also natural to consider the approximation of continuous maps by regular maps. Let $X$, $Y$ be real algebraic varieties (possibly singular), and $D$ an arbitrary subset of $X$. We say that a continuous map $f:D\rightarrow Y$ can be \emph{approximated in the compact-open topology by regular maps from $X$ to $Y$} if for every open neighborhood $\mathcal{U}\subset \mathcal{C}(D,Y)$ of $f$ there exists a regular map $g:X\rightarrow Y$ whose restriction $g|_D$ belongs to $\mathcal{U}$. Here $\mathcal{C}(D,Y)$ denotes the space of all continuous maps from $D$ to $Y$, endowed with the compact-open topology.

\begin{theorem}
\label{thm:1.12}
    Let $Y$ be a quasi-malleable nonsingular real algebraic variety. Let $X$ be a real algebraic variety (possibly singular), and $f:D\rightarrow Y$ a continuous map defined on an arbitrary subset $D$ of $X$. Assume that there exists a regular map $f_0:X\rightarrow Y$ with $f_0|_D$ homotopic to $f$. Then $f$ can be approximated in the compact-open topology by regular maps from $X$ to $Y$.
\end{theorem}

The following is an immediate corollary to Theorem \ref{thm:1.12}.

\begin{corollary}
\label{cor:1.13}
    Let $Y$ be a quasi-malleable nonsingular real algebraic variety. Let $X$ be a real algebraic variety (possibly singular), and $D$ an arbitrary contractible subset of $X$. Then every continuous map from $D$ to $Y$ can be approximated in the compact-open topology by regular maps form $X$ to $Y$.
\end{corollary}

We also have the following variant of Theorem \ref{thm:1.12} that takes interpolation into account.

\begin{theorem}
\label{thm:1.14}
    Let $Y$ be a malleable nonsingular real algebraic variety. Let $X$ be a real algebraic variety (possibly singular), $f:D\rightarrow Y$ a continuous map defined on an arbitrary subset $D$ of $X$, and $A$ a Zariski closed subvariety of $X$. Assume that there exists a regular map $f_0:X\rightarrow Y$ with $f_0|_{D\cap A}=f|_{D\cap A}$ and with $f_0|_D$ homotopic to $f$ relative to $D\cap A$. Then $f$ can be approximated in the compact-open topology by regular maps from $X$ to $Y$ whose restrictions to $A$ are equal to $f_0|_A$.
\end{theorem}

This paper focuses on the AP and AIP of the target varieties. Other aspects of the problem of approximation by regular maps are studied in \cite{5,7,8,9,10,11,12,13,14,15,16,18,26,27,28,32,33,37,43,45,47,49,50,51,52,53,56,58}; the recent paper \cite{2} by Benoist and Wittenberg contains important results for maps defined on real algebraic curves. The literature also considers approximation by more flexible "almost regular" maps, where "almost regular" has a precise context-dependent meaning (see \cite{3,4,36,38,39,41,42,43,49,59} and the references therein). Remarkable results on approximation by differentiable semialgebraic maps are given by Fernando and Ghiloni \cite{23} as well as by Pawłucki \cite{52}.

The paper is organized as follows. Working with smooth manifolds and smooth maps, we develop our main tools in Section \ref{section:2}. These allow us to prove Theorems \ref{thm:1.3} and \ref{thm:1.12} in Section \ref{section:3}. In Section \ref{section:4}, we prove Propositions \ref{prop:4.1} and \ref{prop:4.4}, which correspond to Theorems \ref{thm:1.6} and \ref{thm:1.14} with $Y=\mathbb{R}$. Combining the results of Sections \ref{section:2} and \ref{section:4}, we prove Theorems \ref{thm:1.6} and \ref{thm:1.14} in Section \ref{section:5}. Section \ref{section:6}, which does not depend on the previous sections, contains the proof of Theorem \ref{thm:1.9}. The inspiration for the use of dominating sprays comes from complex geometry, primarily from Gromov's paper \cite{29} and the related works of Forstneri\v{c} and others elaborated in \cite{25}.

\section{Computations in smooth manifolds}
\label{section:2}

The proofs of Theorems \ref{thm:1.3}, \ref{thm:1.6}, \ref{thm:1.12} and \ref{thm:1.14} depend largely on topological arguments. This section is devoted to developing the necessary topological tools.

\begin{notation}
    Given a homotopy $F:X\times [0,1]\rightarrow Y$ of continuous maps of topological spaces and a number $t\in [0,1]$, we denote by $F_t$ the map from $X$ to $Y$ defined by $F_t(x)=F(x,t)$, for all $x\in X$. If $T$ is a subset of a metric space $S$ and $\varepsilon>0$ is a constant, then the \emph{$\varepsilon$-neighborhood} of $T$ in $S$ is the union of all open balls with centers in the points of $T$ and radius $\varepsilon$. By a metric on a smooth manifold we always mean a metric which induces the manifold topology. 
\end{notation} 

For the general theory of smooth vector bundles over smooth manifolds we refer to \cite{30}. The notions of quasi-spray and spray introduced in Definition \ref{def:1.2} have natural variants in the category of smooth manifolds and smooth maps. To be precise, let $N$ be a smooth manifold. A \emph{smooth quasi-spray} $(E,p,\{s_V\},\{\tau_V\})$ for $N$ consists of a smooth vector bundle $p:E\rightarrow N$ over $N$ and two collections $\{s_V\}$, $\{\tau_V\}$ indexed by $V\in \Omega(N)$, where $s_V:E\rightarrow N$ and $\tau_V:V\rightarrow E$ are smooth maps such that $p(\tau_V(y))=y$ and $s_V(\tau_V(y))=y$ for all $y\in V$. A smooth quasi-spray $(E,p,\{s_V\},\{\tau_V\})$ for $N$ is said to be \emph{dominating} if for every $V\in\Omega(N)$ and every $y\in V$ the smooth map 
$$E_y\rightarrow N,\quad v\mapsto s_V(v+\tau_V(y))$$
is a submersion at $0_y\in E_y$. A \emph{smooth spray} $(E,p,s)$ for $N$ consists of a smooth vector bundle $p:E\rightarrow N$ over $N$ and a smooth map $s:E\rightarrow N$ such that $s(0_y)=y$ for all $y\in N$. A smooth spray $(E,p,s)$ for $N$ is said to be \emph{dominating} if for every $y\in N$ the smooth map 
$$E_y\rightarrow N,\quad v\mapsto s(v)$$
is a submersion at $0_y\in E_y$.

We need an additional definition to work effectively with dominating quasi-sprays.

\begin{definition}
\label{def:2.2}
    Let $N$ be a smooth manifold.
    \begin{enumerate}
        \item[(1)] (Smooth partial spray).  A \emph{smooth partial spray} $(E,E^{\circ},p,s)$ for $N$ consists of a smooth vector bundle $p:E\rightarrow N$ over $N$, an open neighborhood $E^{\circ}\subset E$ of the zero section $Z(E)$, and a smooth map $s:E^{\circ}\rightarrow N$ such that $s(0_y)=y \textnormal{ for all }y\in N$.
        \item[(2)] (Dominating smooth partial spray).  A smooth partial spray $(E,E^{\circ},p,s)$ for $N$ is said to be \emph{dominating} if for every point $y\in N$ the smooth map
        $$E_y\cap E^{\circ}\rightarrow N,\quad v\mapsto s(v)$$
        is a submersion at $0_y\in E_y$.
    \end{enumerate}
\end{definition}

A smooth partial spray $(E,E^{\circ},p,s)$ where $E^{\circ}=E$ is identified with the spray $(E,p,s)$.

If $E$ is a smooth vector bundle over a smooth manifold $N$, then for each point $y\in N$ the tangent space $T_{0_y}E$ to $E$ at $0_y\in E_y$ splits as 
$$T_{0_y}E=T_{0_y}Z(E)\oplus E_y$$
where $E_y$ is identified with $T_{0_y}E_y\subset T_{0_y}E$. Moreover, a smooth partial spray $(E,E^{\circ},p,s)$ for $N$ is dominating if and only if for every point $y\in N$ the derivative 
$$d_{0_y}s:T_{0_y}E^{\circ}=T_{0_y}E\rightarrow T_yN$$
maps $E_y$ onto $T_yN$.

\begin{lemma}
\label{lem:2.3}
    Let $N$ be a smooth manifold and let $(E,E^{\circ},p,s)$ be a dominating smooth partial spray for $N$. Then there exists a dominating smooth partial spray $(\tilde{E},\tilde{E^{\circ}},\tilde{p},\tilde{s})$ for $N$, where $\tilde{E}$ is a smooth vector subbundle of $E$, $\tilde{E}^{\circ}=\tilde{E}\cap E^{\circ}$, $\tilde{p}=p|_{\tilde{E}}$ and $\tilde{s}=s|_{\tilde{E}^{\circ}}$, such that $\tilde{E}$ regarded as a smooth manifold has dimension $2\textnormal{dim}N$. Moreover, the map 
    $$\tilde{E}\rightarrow TN,\quad v\mapsto d_{0_{\tilde{p}(v)}}\tilde{s}(v)$$
    is a smooth isomorphism of smooth vector bundles.
\end{lemma}
\begin{proof}
    Since $(E,E^{\circ},p,s)$ is a dominating smooth partial spray for $N$, the map 
    $$E\rightarrow TN,\quad v\mapsto d_{0_{p(v)}}s(v)$$
    is a surjective smooth morphism of smooth vector bundles, and hence its kernel $F$ is a smooth vector subbundle of $E$. Therefore $E$ splits as a direct sum $E=\tilde{E}\oplus F$ for some smooth vector subbundle $\tilde{p}:\tilde{E}\rightarrow N$ of $p:E\rightarrow N$. If $\tilde{E}^{\circ}=\tilde{E}\cap E^{\circ}$ and $\tilde{s}=s|_{\tilde{E}^{\circ}}$, then 
    $$d_{0_{\tilde{p}(v)}}\tilde{s}(v)=d_{0_{p(v)}}s(v)\, \textnormal{for all }v\in \tilde{E},$$
    which completes the proof.
\end{proof}

The following lemma will allow us to avoid having to repeat the same argument a number of times later on.

\begin{lemma}
\label{lem:2.4}
    Let $N$ be a smooth manifold and let $(E,E^{\circ},p,s)$ be a dominating smooth partial spray for $N$. Assume that $E$ regarded as a smooth manifold has dimension $2$\textnormal{dim}$N$. Then the smooth map 
    $$\Phi:=(p,s):E^{\circ}\rightarrow N\times N$$
    induces a diffeomorphism from the zero section $Z(E)$ of $E$ onto the diagonal $\triangle_{N}$ of $N\times N$. Moreover, $\Phi$ maps an open neighborhood of $Z(E)$ in $E^{\circ}$ diffeomorphically onto an open neighborhood of $\triangle_N$ in $N\times N$.
\end{lemma}
\begin{proof}
    The first assertion holds because $\Phi(0_y)=(y,y)$ for all $y\in N$. We claim that for each $y\in N$ the derivative
    $$d_{0_y}\Phi:T_{0_y}E=T_{0_y}Z(E)\oplus E_y\rightarrow T_{(y,y)}(N\times N)=T_yN\times T_yN$$
    is an isomorphism. It suffices to show that it is injective because the tangent spaces $T_{0_y}E$ and $T_{(y,y)}(N\times N)$ have the same dimension. If $d_{0_y}\Phi(w)=0$ for some vector $w\in T_{0_y}E$, then $d_{0_y}p(w)=0$ and $d_{0_y}s(w)=0$. Since the kernel of $d_{0_y}p$ is precisely $E_y$ and $d_{0_y}s$ is injective when restricted to $E_y$, we get $w=0$. Thus, $d_{0_y}\Phi$ is injective as required. Now it is enough to invoke \cite[(12.7)]{20} to complete the proof.
\end{proof}

The following lemma is key.

\begin{lemma}
\label{lem:2.5}
    Let $N$ be a smooth manifold endowed with a metric \textnormal{dist}. Assume that $(E,E^{\circ},p,s)$ is a dominating smooth partial spray for $N$. Let $M$ be a smooth manifold, $M_0\subset M$ an open set (possibly empty), $B\subset M$ an arbitrary set and $K\subset M$ a compact set containing both $M_0$ and $B$. Let $f:K\rightarrow N$ be a continuous map whose restriction $f|_{M_0}$ is smooth. Then there exists a constant $\varepsilon>0$ such that for every continuous map $g:K\rightarrow N$ satisfying 
    \begin{named}
        \item \textnormal{dist}$(f(x),g(x))<\varepsilon \,$ for all $x\in K$,
        \item $g|_{M_0}$ is smooth,
        \item $g|_B=f|_B$,
    \end{named}
    there exists a continuous map $\xi:K\rightarrow E^{\circ}$ satisfying
    \begin{named}
        \item $p\circ \xi =f$,
        \item $s\circ \xi =g$,
        \item $\xi|_{M_0}$ is smooth,
        \item $\xi(x)=0_{f(x)}$ for all $x\in B$.
    \end{named}
\end{lemma}
\begin{proof}
    By Lemma \ref{lem:2.3}, we may assume that the map 
    $$E\rightarrow TN,\quad v\mapsto d_{0_{p(v)}}s(v)$$
    is a smooth isomorphism of smooth vector bundles. In particular, $E$ regarded as a smooth manifold has dimension $2$dim$N$.

    Consider the smooth map
    $$\Phi :=(p,s):E^{\circ}\rightarrow N\times N.$$
    According to Lemma \ref{lem:2.4}, there exist an open neighborhood $U\subset E^{\circ}$ of $Z(E)$ and an open neighborhood $V\subset N\times N$ of the diagonal $\triangle_N$ such that $\Phi(U)=V$ and $\Phi|_{U}:U\rightarrow V$ is a diffeomerphism. Since the continuous map $(f,f):K\rightarrow N\times N$ transforms $K$ onto a compact set $(f,f)(K)$ in $\triangle_N$, there exists a constant $\varepsilon>0$ such that the $\varepsilon$-neighborhood of $(f,f)(K)$ in $N\times N$ is contained in $V$; here we refer to the $\varepsilon$-neighborhood with respect to the metric $d$ on $N\times N$ defined by 
    $$d((y_1,z_1),(y_2,z_2))=\textnormal{max}\{ \textnormal{dist}(y_1,y_2),\textnormal{dist}(z_1,z_2) \}.$$
    Now let $g:K\rightarrow N$ be a continuous map satisfying conditions (2.5.1)-(2.5.3). By the definition of $\varepsilon$, the set $(f,g)(K)$ is contained in $V$, so we can define $\xi :K\rightarrow E^{\circ}$ as $\xi=(\Phi|_{U})^{-1}\circ (f,g)$. Conditions (2.5.4)-(2.5.7) follow immediately.
\end{proof}

In the proof of Theorem \ref{thm:1.3}, the following variant of Lemma \ref{lem:2.5} will play an important role.

\begin{lemma}
\label{lem:2.6}
    Let $N$ be a smooth manifold endowed with a metric \textnormal{dist}. Assume that $(E,p,\{s_V\},\{\tau_V\})$ is a dominating quasi-spray for $N$. Let $M$ be a smooth manifold, $M_0\subset M$ an open set (possibly empty) and $K\subset M$ a compact set containing $M_0$. Let $f:K\rightarrow N$ be a continuous map whose restriction $f|_{M_0}$ is smooth. Then there exist a constant $\varepsilon>0$ and a set $W\in\Omega(N)$ such that for every continuous map $g:K\rightarrow N$ satisfying 
    \begin{named}
        \item \textnormal{dist}$(f(x),g(x))<\varepsilon$ for all $x\in K$,
        \item $g|_{M_0}$ is smooth,
    \end{named}
    there exists a continuous map $\xi:K\rightarrow E$ satisfying
    \begin{named}
        \item $p\circ \xi = f$,
        \item $s_W\circ \xi = g$,
        \item $\xi|_{M_0}$ is smooth.
    \end{named}
\end{lemma}
\begin{proof}
    Since $f(K)$ is a compact subset of $N$, there exists an open set $W\in \Omega(N)$ that contains $f(K)$. We choose a constant $\delta>0$ such that the $\delta$-neighborhood of $f(K)$ in $N$ is contained in $W$. Thus, if $g:K\rightarrow N$ is a map satisfying 
    $$d(f(x),g(x))<\delta\quad \textnormal{for all }x\in K$$
    then $g(K)\subset W$.

    Let $E(W):=p^{-1}(W)$ be the restriction of the vector bundle $E$ to $W$. For each vector $v\in E(W)$, both vectors $v$ and $\tau_W(p(v))$ belong to $E_{p(v)}$, so the map
    $$\hat{s}_W:E(W)\rightarrow N,\quad v\mapsto s_W(v+\tau_W(p(v)))$$
    is well defined and smooth, with $\hat{s}_W(0_{p(v)})=p(v)$. Thus, defining 
    $$E(W)^{\circ}:=(\hat{s}_W)^{-1}(W)$$
    $$s_W^{\circ}:E(W)^{\circ}\rightarrow W\subset N,\quad s_W^{\circ}:=\hat{s}_W|_{E(W)^{\circ}}$$
    we obtain a dominating smooth partial spray $(E(W),E(W)^{\circ},p|_{E(W)},s_W^{\circ})$ for $W$. By Lemma \ref{lem:2.5} (with $B=\varnothing$), there exists a constant $\varepsilon$, $0<\varepsilon\leq \delta$, such that if $g:K\rightarrow N$ is a continuous map satisfying conditions (2.6.1) and (2.6.2), then there exists a continuous map $\xi_W:K\rightarrow E(W)^{\circ}$ with the following properties:
    \begin{itemize}
        \item $p\circ \xi_W=f$,
        \item $s_W^\circ\circ \xi_W=g$,
        \item $\xi_W|_{M_0}$ is smooth.
    \end{itemize}
    Now it is enough to define $\xi$ as $\xi=\xi_W+\tau_W\circ f$. Conditions (2.6.3)-(2.6.5) are easy to check.
\end{proof}

We now turn to auxiliary results of a different kind, which are needed for the proofs of our main theorems.

Let $\zeta$ be a smooth vector field on a smooth manifold $N$. Since $\zeta$ induces a smooth local flow, we can choose an open neighborhood $D\subset N\times \mathbb{R}$ of $N\times \{ 0 \}$ and a smooth map $\theta:D\rightarrow N$ satisfying
$$\theta(y,0)=y \textnormal{ and } \frac{\partial\theta}{\partial t}(y,0)=\zeta(y) \textnormal{ for all } y\in N,$$
where $\frac{\partial\theta}{\partial t}(y,0)$ is the tangent vector to the parametric curve $t\mapsto \theta(y,t)$ at $t=0$.

A slightly stronger statement is also true.

\begin{lemma}
\label{lem:2.7}
    If $\zeta$ is a smooth vector field on a smooth manifold $N$, then there exists a smooth map $\phi:N\times\mathbb{R}\rightarrow N$ such that 
    $$\varphi(y,0)=y \textnormal{ and } \frac{\partial\varphi}{\partial t}(y,0)=\zeta(y) \textnormal{ for all } y\in N.$$
\end{lemma}
\begin{proof}
    Let $\theta:D\rightarrow N$ be as above. Choose a smooth function $\varepsilon:N\rightarrow (0,\infty)$ such that the set $\{ y \}\times (-\varepsilon(y),\varepsilon(y))$ is contained in $D$ for all $y\in N$. Then the map $\varphi:N\times\mathbb{R}\rightarrow N$ defined by
    $$\varphi(y,t)=\theta(y,\frac{\varepsilon(y)t}{\sqrt{\varepsilon(y)^2+t^2}}) \textnormal{ for all } (y,t)\in N\times \mathbb{R}$$
    has the required properties.
\end{proof}

We say that $n$ smooth vector fields $\zeta_1$, \ldots, $\zeta_n$ on a smooth manifold $N$ \emph{generate} the tangent bundle $TN$ if for every point $y\in N$ the vectors $\zeta_1(y)$, \ldots, $\zeta_n(y)$ span the tangent space $T_yN$. If $N$ is a smooth submanifold of $\mathbb{R}^n$, then $TN$ is generated by $n$ smooth vector fields because there exists a smooth vector bundle morphism from the product vector bundle $N\times \mathbb{R}^n$ onto $TN$.

\begin{lemma}
\label{lem:2.8}
    Let $N$ be a smooth manifold whose tangent bundle $TN$ is generated by $n$ smooth vector fields. Then there exist $n$ smooth maps $\varphi^1$, \ldots, $\varphi^n$ from $N\times \mathbb{R}$ to $N$ such that 
    $$\varphi^i(y,0)=y \textnormal{ for all } i\in \{ 1, \ldots, n \} \textnormal{ and all } y\in N,$$
    and the vector fields 
    $$\frac{\partial\varphi^1}{\partial t}(\cdot,0), \ldots, \frac{\partial\varphi^n}{\partial t}(\cdot,0)$$
    on $N$ generate $TN$. Moreover, $(E,p,\sigma)$ is a dominating smooth spray for $N$, where $p:E=N\times \mathbb{R}^n\rightarrow N$ is the product vector bundle and $\sigma:E\rightarrow N$ is defined by 
    $$\sigma(y,(t_1,\ldots,t_n))=(\varphi_{t_n}^n\circ\ldots\circ\varphi_{t_1}^1)(y) \textnormal{ for all } (y,(t_1,\ldots,t_n))\in N\times\mathbb{R}^n$$
    with $\varphi_{t_i}^i=\varphi^i(\cdot,t_i)$ for all $i\in \{ 1, \ldots, n \}$.
\end{lemma}
\begin{proof}
    The first assertion follows from Lemma \ref{lem:2.7}. A straightforward computation shows that the second assertion holds.
\end{proof}

For the sake of completeness we include the following observation.

\begin{lemma}
\label{lem:2.9}
    Let $f:M\rightarrow N$ be a smooth embedding of smooth manifolds, where \textnormal{dim}$M$=\textnormal{dim}$N$. Let $C\subset M$ be a compact set, and $M_0\subset M$ a relatively compact open set containing $C$. Then $f$ has an open neighborhood $\mathcal{U}\subset \mathcal{C}^{\infty}(M,N)$ in the $\mathcal{C}^1$ topology such that for every map $g\in \mathcal{U}$ the restriction $g|_{M_0}$ is a smooth embedding and $f(C)\subset g(M_0)$.
\end{lemma}
\begin{proof}
    Choose a relatively compact open set $U\subset M$ containing the closure $\bar{M_0}$ of $M_0$. Let $\varphi_0$, $\varphi_1$ be two smooth real-valued functions on $M$ satisfying 
    $$M\setminus U\subset \varphi^{-1}_0(0),\, \bar{M_0}\subset \varphi^{-1}_0(1) \textnormal{ and } M\setminus M_0\subset \varphi^{-1}_1(0),\, C\subset \varphi^{-1}_1(1).$$
    By Sard's theorem there exists a regular value $\lambda_i\in \mathbb{R}$ for $\varphi_i$ with $0<\lambda_i<1$ for $i=1,2$. Then 
    $$K_i:=\left\{ x\in M:\quad \varphi_i(x)\geq \lambda_i \right\}\subset M$$
    is a compact smooth submanifold with boundary; the boundary $\partial K_i=\varphi^{-1}_i(\lambda_i)$ of $K_i$ can be empty. Moreover,
    $$C\subset K_1\setminus\partial K_1\subset K_1\subset M_0\subset \bar{M_0}\subset K_0\subset M.$$
    By \cite[p.37, Theorem 1.4]{30}, $f$ has an open neighborhood $\mathcal{U}\subset \mathcal{C}^{\infty}(M,N)$ in the $\mathcal{C}^1$ topology such that for every map $g\in \mathcal{U}$ the restriction $g|_{K_0}$ is a smooth embedding, so $g|_{M_0}$ is also a smooth embedding. Shrinking $\mathcal{U}$ if necessary, we get $f(C)\subset g(K_1\setminus \partial K_1)\subset g(M_0)$ in view of Lemma \ref{lem:2.10} below.
\end{proof}

In the proof of Lemma \ref{lem:2.9} we used the following lemma.

\begin{lemma}
\label{lem:2.10}
    Let $K$ be a compact smooth manifold with boundary, $C$ a compact set in $K\setminus\partial K$, and $N$ a smooth manifold. Let $f:K\rightarrow N$ be a smooth embedding. Assume that \textnormal{dim}$K=$\textnormal{dim}$N$. Then $f$ has an open neighborhood $\mathcal{U}\subset \mathcal{C}^{\infty}(K,N)$ in the $\mathcal{C}^1$ topology such that every map $g\in \mathcal{U}$ is a smooth embedding satisfying $f(C)\subset g(K\setminus \partial K)$.
\end{lemma}
\begin{proof}
    Choose two disjoint open sets $V$ and $W$ in $N$ such that $f(C)\subset V$ and $f(\partial K)\subset W$. For each point $x\in C$ choose a path connected open set $V_x$ in $N$ with $f(x)\in V_x\subset V$. Let $U_x\subset N$ be an open neighborhood of $f(x)$ whose closure $\bar{U}_x$ is compact and contained in $V_x$. The compact set $C$ is covered by the open sets $f^{-1}(U_x)$, where $x\in C$, and hence there exist points $x_1$, \ldots, $x_r$ in $C$ such that $C$ is also covered by the sets $f^{-1}(U_{x_i})$, where $i\in \{ 
    1, \ldots, r \}$. Let $U_i:=U_{x_i}$, $V_i:=V_{x_i}$ and $C_i:=C\cap f^{-1}(\bar{U_i})$. Obviously, the $C_i$ are compact sets covering $C$.

    By \cite[p.37, Theorem 1.4]{30}, $f$ has an open neighborhood $\mathcal{U}\subset \mathcal{C}^{\infty}(K,N)$ in the $\mathcal{C}^1$ topology such that every map $g\in \mathcal{U}$ is a smooth embedding. Shrinking $\mathcal{U}$ if necessary, we get additionally that 
    $$g(\partial K)\subset W \textnormal{ and } g(C_i)\subset V_i \textnormal{ for } i\in \{ 1,\ldots, r \}.$$
    Now, for a given $g\in \mathcal{U}$, it is sufficient to check that the inclusion $f(C)\subset g(K\setminus \partial K)$ holds. If $x\in C$, then $x\in C_i$ for some $i\in \{ 1,\ldots r \}$, and hence both $f(x)$ and $g(x)$ belong to $V_i$. Therefore there exists a continuous path $\gamma:[0,1]\rightarrow V_i$ with $\gamma(0)=g(x)$ and $\gamma(1)=f(x)$. Since $K$ is a compact space, its image $g(K)$ is closed in $N$. On the other hand, $g$ is a smooth embedding, so $g(K)$ is the union of the open set $g(K\setminus \partial K)$ and the set $g(\partial K)$ disjoint from $V$. It follows that the subset 
    $$I:=\gamma^{-1}(g(K))=\gamma^{-1}(g(K\setminus \partial K))$$
    of the interval $[0,1]$ is both open and closed. Since $0\in I$, we get $I=[0,1]$ due to the connectedness of $[0,1]$. So $f(x)\in g(K\setminus\partial K)$, which completes the proof.
\end{proof}

\section{Real algebraic varieties with the AP}
\label{section:3}

This section is devoted to the proofs of Theorems \ref{thm:1.3} and \ref{thm:1.12}.

\begin{lemma}
\label{lem:3.1}
    Let $Y$ be a quasi-malleable nonsingular real algebraic variety. Then there exists a dominating quasi-spray $(E,p,\{s_V\},\{\tau_V\})$ for $Y$ such that $p:E=Y\times \mathbb{R}^n\rightarrow Y$ is the product vector bundle for some $n$.
\end{lemma}
\begin{proof}
    Let $(E',p',\{s_V'\},\{\tau_V'\})$ be a dominating quasi-spray for $Y$. Choose an algebraic vector bundle $E''$ over $Y$ such that the direct sum $E'\oplus E''$ is algebraically isomorphic to the product vector bundle $p:Y\times \mathbb{R}^n\rightarrow Y$ for some nonnegative integer $n$ (see \cite[Theorem 12.17]{5}). Identifying $E'\oplus E''$ with $Y\times \mathbb{R}^n$, we define $s_V:Y\times \mathbb{R}^n\rightarrow Y$ to be the composite of the projection $E'\oplus E''\rightarrow E'$ with $s_V':E'\rightarrow Y$, and $\tau_V:V\rightarrow Y\times \mathbb{R}^n$ to be the composite of $\tau_V':V\rightarrow E'$ with the inclusion $E'\hookrightarrow E'\oplus E''$. By construction, $(E,p,\{s_V\},\{\tau_V\})$ is a dominating quasi-spray for $Y$.
\end{proof}

\begin{theorem}
\label{thm:3.2}
    Let $Y$ be a quasi-malleable nonsingular real algebraic variety. Then $Y$ has the generalized approximation property.
\end{theorem}
\begin{proof}
    In the proof, all statements about approximation refer to approximation in the $\mathcal{C}^{\infty}$ topology.

    Let $X$ be a nonsingular real algebraic variety, $U$ an open subset of $X$, and $f:U\rightarrow Y$ a smooth map. Assume that there exists a regular map $f_0:X\rightarrow Y$ with $f_0|_U$ homotopic to $f$. Our goal is to prove that $f$ can be approximated by regular maps from $X$ to $Y$. In order to achieve this, it is sufficient to show that for any given open set $U_0\in\Omega(U)$, the restriction $f|_{U_0}$ can be approximated by regular maps from $X$ to $Y$.

    Choose a homotopy $F:U\times[0,1]\rightarrow Y$ with $F_0=f_0|_U$ and $F_1=f$. By \cite[Proposition 10.22]{44}, we may assume that $F$ is a smooth map. Consider the subset 
    $$I:=\{ t\in[0,1]:\quad F_t|_{U_0} \textnormal{ can be approximated by regular maps from } X \textnormal{ to } Y \}$$
    of $[0,1]$. By definition, $I$ is closed and $0\in I$. We now prove that $I$ is open in $[0,1]$. Let $t_0\in I$. Let dist be a metric on $Y$. According to Lemma \ref{lem:3.1}, there exists a dominating quasi-spray $(E,p,\{s_V\},\{\tau_V\})$ for $Y$, where $p:E=Y\times \mathbb{R}^n\rightarrow Y$ is the product vector bundle over $Y$. Choose a constant $\varepsilon>0$ and an open set $W\in\Omega(Y)$ as in Lemma \ref{lem:2.6} (with $N=Y$, $M=U$, $M_0=U_0$, $K=\bar{U_0}$, $f=F_{t_0}|_K$). Now choose a neighborhood $I_0$ of $t_0$ in $[0,1]$ such that 
    $$\textnormal{dist}(F_t(x),F_{t_0}(x))<\varepsilon \textnormal{ for all } t\in I_0 \textnormal{ and all } x\in K.$$
    Thus, for each $t\in I_0$ there exists a continuous map $\xi:K\rightarrow E=Y\times \mathbb{R}^n$ satisfying
    \begin{enumerate}
        \item[(i)] $p\circ\xi = F_{t_0}|_K$,
        \item[(ii)] $s_W\circ\xi=F_t|_K$,
        \item[(iii)] $\xi|_{U_0}$ is smooth.
    \end{enumerate}
    By (i), the map $\xi$ is of the form $\xi=(F_{t_0}|_K,\eta)$ for some continuous map $\eta:K\rightarrow\mathbb{R}^n$. In view of (iii), the restriction $\eta|_{U_0}$ is a smooth map. Moreover, using (ii) we get 
    $$F_t(x)=s_W(F_{t_0}(x),\eta(x)) \textnormal{ for all } x\in U_0.$$
    According to the definition of $I$, the map $F_{t_0}|_{U_0}$ can be approximated by regular maps from $X$ to $Y$. By the Stone-Weierstrass theorem, the map $\eta$ can be approximated by regular maps from $X$ to $\mathbb{R}^n$. Thus, since $s_W$ is a regular map, we obtain that $F_t|_{U_0}$ can be approximated by regular maps from $X$ to $Y$, so $t\in I$. Therefore the inclusion $I_0\subset I$ holds, which means that $I$ is an open subset of $[0,1]$. In summary, $I=[0,1]$ due to the connectedness of the interval $[0,1]$. The proof is complete because $f|_{U_0}=F_1|_{U_0}$.
\end{proof}
\begin{theorem}
\label{thm:3.3}
    Let $Y$ be a nonsingular real algebraic variety. Let $(E,p,\sigma)$ be a dominating smooth spray for $Y$, where $p:E\rightarrow Y$ is an algebraic vector bundle over $Y$. Assume that the smooth map $\sigma:E\rightarrow Y$ can be approximated in the $\mathcal{C}^1$ topology by regular maps from $E$ to $Y$. Then $Y$ admits a dominating quasi-spray of the form $(E,p,\{s_V\},\{\tau_V\})$. In particular, the variety $Y$ is quasi-malleable.
\end{theorem}
\begin{proof}
    By Lemma \ref{lem:2.3} (with $N=Y$, $E^{\circ}=E$), there exists a dominating smooth spray $(\tilde{E},\tilde{p},\tilde{\sigma})$ for $Y$ such that $\tilde{p}:\tilde{E}\rightarrow Y$ is a smooth vector subbundle of the algebraic vector bundle $p:E\rightarrow Y$, $\tilde{\sigma}=\sigma|_{\tilde{E}}$, and $\tilde{E}$ regarded as a smooth manifold has dimension $2$dim$Y$. Hence, according to Lemma \ref{lem:2.4}, the smooth map 
    $$\Phi:=(\tilde{p},\tilde{\sigma})=(p|_{\tilde{E}},\sigma|_{\tilde{E}}):\tilde{E}\rightarrow Y\times Y$$
    maps an open neighborhood $M\subset \tilde{E}$ of the zero section $Z(\tilde{E})$ diffeomorphically onto an open neighborhood $\Phi(M)\subset Y\times Y$ of the diagonal $\triangle_Y$ of $Y\times Y$.

    Let $V\in \Omega(Y)$. The closure $\bar{V}$ of $V$ is a compact set in $Y$, so the diagonal $\triangle_{\bar{V}}$ of $\bar{V}\times \bar{V}$ is a compact set in $\triangle_Y$, and the preimage $C:=(\Phi|_M)^{-1}(\triangle_{\bar{V}})$ is a compact set in $M$. Choose a relatively compact open set $M_0\subset M$ that contains $C$. 

    Let $s_V:E\rightarrow Y$ be a regular map. We consider the smooth map 
    $$\Psi:=(p|_M,s_V|_M):M\rightarrow Y\times Y.$$
    If $s_V$ is sufficiently close to $\sigma$ in the $\mathcal{C}^1$ topology, then by Lemma \ref{lem:2.9} the restriction $\Psi|_{M_0}$ is a smooth embedding such that $\Psi(C)=\triangle_{\bar{V}}\subset \Psi(M_0)$. Denoting by $\psi:M_0\rightarrow \Psi(M_0)$ the diffeomorphism induced by $\Psi$, we define a map $\tau:V\rightarrow \bar{E}$ by $\tau(y)=\psi^{-1}(y,y)$ for all $y\in V$. By construction, $\tau$ is a smooth map satisfying $p(\tau(y))=y$ and $s_V(\tau(y))=y$ for all $y\in V$. Moreover, for each $y\in V$, the derivative of the regular map
    $$\tilde{E}_y\rightarrow Y,\quad v\mapsto s_V(v+\tau(y))$$
    at $0_y\in \tilde{E}_y$ is an isomorphism since $\psi$ is a diffeomorphism and 
    $$\psi(v+\tau(y))=(p(v+\tau(y)),s_V(v+\tau(y)))=(y,s_V(v+\tau(y)))$$
    for all $v$ near $0_y\in\tilde{E}_y$. Define $\tau_V:V\rightarrow E$ to be the composite of $\tau:V\rightarrow \tilde{E}$ and the inclusion vector bundle morphism $\tilde{E}\hookrightarrow E$. Now, if we let $V$ vary in $\Omega(Y)$, we get a dominating quasi-spray for $Y$, so $Y$ is a quasi-malleable variety.
\end{proof}

The following observation is included for the sake of clarity.

\begin{observation}
\label{obs:3.4}
    Let $p:E\rightarrow Y$ be an algebraic vector bundle over a real algebraic variety $Y$ and let $f:E\rightarrow Y$ be a continuous map such that $f(0_y)=y$ for all $y\in Y$. Then $f$ is homotopic to a regular map from $E$ to $Y$. Actually, 
    $$F:E\times [0,1]\rightarrow Y,\quad F(v,t)=f(tv)$$
    is a homotopy from $p$ to $f$ relative to the zero section $Z(E)$.
\end{observation}

We are now in position to carry out the main tasks of this section.

\begin{proof}[Proof of Theorem \ref{thm:1.3}]
    Obviously, (a) implies (b). In light of Observation \ref{obs:3.4}, (b) implies both (c) and (d). By Theorem \ref{thm:3.2}, (e) implies (a).

    According to Lemma \ref{lem:2.8} and Lemma \ref{lem:2.3} (with $E^{\circ}=E$), there exists a dominating smooth spray for $Y$ of the form $(TY,p,\sigma)$, where $p:TY\rightarrow Y$ is the tangent bundle of $Y$. Therefore (c) implies (e) by Theorem \ref{thm:3.3}.

    Suppose (d) holds. Then, by Lemma \ref{lem:2.8}, there exists a dominating smooth spray $(E,p,\sigma)$ for $Y$, where $p:E=Y\times \mathbb{R}^n\rightarrow Y$ is the product vector bundle and the smooth map $\sigma:E\rightarrow Y$ can be approximated in the $\mathcal{C}^1$ topology by regular maps from $E$ to $Y$. Hence (e) holds by Theorem \ref{thm:3.3}.
\end{proof}

It will be seen that the proof of Theorem \ref{thm:1.12} given below is analogous to the proof of Theorem \ref{thm:3.2}.

\begin{proof}[Proof of Theorem \ref{thm:1.12}]
    In the proof, all statements about approximation refer to approximation in the compact-open topology. Our goal is to show that for any given compact subset $K$ of $D$, the restriction $f|_K$ can be approximated by regular maps from $X$ to $Y$. We may assume that $X$ is a Zariski closed subvariety of $\mathbb{R}^m$ for some $m$.

    Let $F:D\times [0,1]\rightarrow Y$ be a homotopy with $F_0=f_0|_D$ and $F_1=f$. Then the subset 
    $$I:=\{ t\in[0,1]:\quad F_t|_K \textnormal{ can be approximated by regular maps from } X \textnormal{ to } Y \}$$
    of $[0,1]$ is closed and $0\in I$. We now prove that $I$ is open in $[0,1]$. Let $t_0\in I$. Let dist be a metric on $Y$. According to Lemma \ref{lem:3.1}, there exists a dominating quasi-spray $(E,p,\{s_V\},\{\tau_V\})$ for $Y$, where $p:E=Y\times \mathbb{R}^n\rightarrow Y$ is the product vector bundle over $Y$. Choose a constant $\varepsilon>0$ and an open set $W\in\Omega(Y)$ as in Lemma \ref{lem:2.6} (with $N=Y$, $M=\mathbb{R}^m$, $M_0=\varnothing$, $f=F_{t_0}|_K$). Now choose a neighborhood $I_0$ of $t_0$ in $[0,1]$ such that 
    $$\textnormal{dist}(F_t(x),F_{t_0}(x))<\varepsilon\quad \textnormal{for all } t\in I_0 \textnormal{ and all } x\in K.$$
    Thus, for each $t\in I_0$ there exists a continuous map $\xi:K\rightarrow E=Y\times \mathbb{R}^n$ satisfying 
    \begin{enumerate}
        \item[(i)] $p\circ\xi=F_{t_0}|_K$,
        \item[(ii)] $s_W\circ\xi=F_t|_K$.
    \end{enumerate}
    By (i), the map $\xi$ is of the form $\xi=(F_{t_0}|_K,\eta)$ for some continuous map $\eta:K\rightarrow \mathbb{R}^n$. Moreover, using (ii), we get 
    $$F_t(x)=s_W(F_{t_0}(x),\eta(x))\quad \textnormal{for all } x\in K.$$
    According to the difinition of $I$, the map $F_{t_0}|_K$ can be approximated by regular maps from $X$ to $Y$. By the Stone-Weierstrass theorem, the map $\eta$ can be approximated by regular maps from $X$ to $\mathbb{R}^n$. Thus, since $s_W$ is a regular map, we obtain that $F_t|_K$ can be approximated by regular maps from $X$ to $Y$, so $t\in I$. Therefore the inclusion $I_0\subset I$ holds, which means that $I$ is an open subset of $[0,1]$. In summary, $I=[0,1]$ due to the connectedness of the interval $[0,1]$. The proof is complete because $f|_K=F_1|_K$.
\end{proof}

\section{Approximation and interpolation of functions}
\label{section:4}

Every real algebraic variety $X$ considered in this paper is affine, and hence the ring $\mathcal{R}(X)$ of all regular functions on $X$ is Noetherian, and furthermore every Zariski closed subvariety $A$ of $X$ is the zero set of the ideal $I_{\mathcal{R}(X)}(A)$ consisting of all functions in $\mathcal{R}(X)$ which vanish on $A$ (see \cite[p.62]{5} or \cite[p.14]{49}).

Assume that $X$ is a nonsingular real algebraic variety. A Zariski closed subvariety $A$ of $X$ is said to be \emph{faithful at a point} $x\in A$ if 
$$I_{\mathcal{O}_x}(A_x)=I_{\mathcal{R}(X)}(A)\mathcal{O}_x,$$
where $\mathcal{O}_x$ is the ring of all real analytic function-germs $(X,x)\rightarrow \mathbb{R}$ and $I_{\mathcal{O}_x}(A_x)$ is the ideal in $\mathcal{O}_x$ comprised of all function-germs vanishing on the germ $A_x$ of $A$ at $x$, while $I_{\mathcal{R}(X)}(A)\mathcal{O}_x$ is the extension of the ideal $I_{\mathcal{R}(X)}(A)$ in $\mathcal{R}(X)$ via the canonical ring homomorphism $\mathcal{R}(X)\rightarrow \mathcal{O}_x$. It readily follows that $A$ is faithful at each of its nonsingular points. The subvariety $A$ is said to be \emph{faithful} if it is faithful at each of its points. Here "faithful" corresponds to "quasi regolare" used in \cite{55}.

First we prove a variant of Tognoli's result \cite[Teorema 1]{22}, generalizing it and omitting the superfluous coherence assumption.

\begin{proposition}
\label{prop:4.1}
    Let $X$ be a nonsingular real algebraic variety, $f:U\rightarrow\mathbb{R}$ a smooth function defined on an open subset $U$ of $X$, and $A$ a Zariski closed subvariety of $X$ which is faithful at every point in $U\cap A$. Assume that there exists a regular function $f_0:X\rightarrow \mathbb{R}$ with $f_0|_{U\cap A}=f|_{U\cap A}$. Then $f$ can be approximated in the $\mathcal{C}^{\infty}$ topology by regular functions on $X$ whose restrictions to $A$ are equal to $f_0|_A$.
\end{proposition}
\begin{proof}
    The function $\varphi:=f-f_0|_U$ is smooth on $U$ and vanishes on $U\cap A$. It suffices to prove that $\varphi$ can be approximated in the $\mathcal{C}^{\infty}$ topology by regular functions on $X$ vanishing on $A$. We may assume that $X$ is a Zariski closed subvariety of $\mathbb{R}^n$ for some $n$. Choose regular functions $h_1,\ldots, h_r$ in $\mathcal{R}(X)$ that generate the ideal $I_{\mathcal{R}(X)}(A)$. Since the subvariety $A$ is faithful at each point in $U\cap A$, we have 
    $$I_{\mathcal{O}_x}(A_x)=I_{\mathcal{R}(X)}(A)\mathcal{O}_x=(h_1,\ldots,h_r)\mathcal{O}_x\quad \textnormal{for all } x\in U\cap A.$$
    It follows that the real analytic subset $U\cap A$ of $U$ is coherent. Thus, by Malgrange's theorem \cite[Chapter VI, Theorem 3.10]{46},
    $$I_{\mathcal{C}_x^{\infty}}(A_x)=I_{\mathcal{O}_x}(A_x)\mathcal{C}^{\infty}_x=(h_1,\ldots h_r)\mathcal{C}_x^{\infty}\quad \textnormal{for all } x\in U\cap A,$$
    where $\mathcal{C}_x^{\infty}$ is the ring of all smooth function-germs $(X,x)\rightarrow \mathbb{R}$, and $I_{\mathcal{C}_x^{\infty}}(A_x)$ is the ideal in $\mathcal{C}_x^{\infty}$ comprised of all function-germs vanishing on $A_x$. Since the function $\varphi$ vanishes on $U\cap A$, its germ $\varphi_x$ at $x$ belongs to the ideal $(h_1,\ldots, h_r)\mathcal{C}_x^{\infty}$ for all $x\in U\cap A$. Hence, using partition of unity, we can express $\varphi$ as
    $$\varphi=a_1h_1|_U+\ldots +a_rh_r|_U$$
    for some smooth functions $a_1,\ldots, a_r$ on $U$. By the Stone-Weierstrass approximation theorem, for each $i\in\{ 
    1,\ldots, r\}$ there exists a regular function $b_i:X\rightarrow\mathbb{R}$ such that $b_i|_U$ is as close to $a_i$ in the $\mathcal{C}^{\infty}$ topology as we want. The function
    $$\psi:=b_1h_1+\ldots +b_rh_r$$
    is regular on $X$ and vanishes on $A$, and $\psi|_U$ is close to $\varphi$ in the $\mathcal{C}^{\infty}$ topology.
\end{proof}

Note that Proposition \ref{prop:4.1} is stronger than the implication (e) $\Rightarrow$ (a) in Theorem \ref{thm:1.6} with $Y=\mathbb{R}$. This is because, using the notation as in Proposition \ref{prop:4.1}, the function $f_0|_U$ is homotopic to $f$ relative to $U\cap A$ via the linear homotopy
$$U\times [0,1]\rightarrow \mathbb{R},\quad (x,t)\mapsto (1-t)f_0(x)+tf(x)$$
and $\mathbb{R}$ is a malleable variety.  

In Proposition \ref{prop:4.1}, the assumption that $A$ is faithful at every point of $U\cap A$ cannot be omitted.

\begin{example}
\label{ex:4.2}
    The real polynomial 
    $$P(x,y)=x^4-2x^2y-y^3$$
    is irreducible and its gradient vanishes only at $(0,0)$. It follows that the real algebraic curve 
    $$C:=\{ (x,y)\in\mathbb{R}^2:\quad P(x,y)=0 \}$$
    is irreducible and has the only singular point at $(0,0)$. However, $C$ is a real analytic submanifold of $\mathbb{R}^2$ because for $y>-1$ the curve is given by the real analytic equation
    $$\varphi(x,y):=x^2-y(1+\sqrt{1+y})=0,$$
    where $\frac{\partial \varphi}{\partial y}(0,0)\neq 0$. Consequently, $C$ is faithful at each point in $C\setminus \{ (0,0) \}$, but is not faithful at $(0,0)$.

    Choose a smooth function $\psi:\mathbb{R}^2\rightarrow \mathbb{R}$ such that $\psi(0,0)\neq 0$ and the support of $\psi$ is contained in the open unit disc in $\mathbb{R}^2$. Then $f:=\psi\varphi$ is a smooth function on $\mathbb{R}^2$ that vanishes on $C$ and satisfies $\frac{\partial f}{\partial y}(0,0)\neq 0$. We claim that $f$ cannot be approximated in the $\mathcal{C}^1$ topology by regular functions vanishing on $C$. Indeed, otherwise we could choose a regular function $g:\mathbb{R}^2\rightarrow\mathbb{R}$ that vanishes on $C$ and satisfies $\frac{\partial g}{\partial y}(0,0)\neq 0$. This would imply that $C$ is nonsingular at $(0,0)$, a contradiction.
\end{example}

As noted in the proof of Proposition \ref{prop:4.1}, the assumption that $A$ is faithful at every point in $U\cap A$ implies the coherence of the real analytic subset $U\cap A$ of $U$. However, if $A$ is faithful at some point $x\in A$, then it is not necessarily coherent at $x$.

\begin{example}
\label{ex:4.3}
    The real algebraic surface 
    $$S:=\{ (x,y,z)\in\mathbb{R}^3:\quad z(x^2+y^2)-x^3=0 \},$$
    called the \emph{Cartan umbrella}, is faithful at $(0,0,0)$, but $S$ regarded as as real analytic set is not coherent at $(0,0,0)$. Note that $S$ is not faithful at any point of the form $(0,0,z)$ with $z\neq 0$.
\end{example}

We now move on to the next result of this section.

\begin{proposition}
\label{prop:4.4}
    Let $X$ be a real algebraic variety (possibly singular), $f:D\rightarrow\mathbb{R}$ a continuous function defined on a compact subset $D$ of $X$, and $A$ a Zariski closed subvariety of $X$. Assume that there exists a regular function $f_0:X\rightarrow\mathbb{R}$ with $f_0|_{D\cap A}=f|_{D\cap A}$. Then $f$ can be approximated in the compact-open topology by regular functions on $X$ whose restrictions to $A$ are equal to $f_0|_A$.
\end{proposition}
\begin{proof}
    The function $\varphi:=f-f_0|_D$ is continuous on $D$ and vanishes on $D\cap A$. It suffices to prove that $\varphi$ can be approximated in the compact-open topology by regular functions on $X$ vanishing on $A$. This amounts to showing that given a compact subset $K$ of $D$ and a constant $\varepsilon>0$, we can choose a regular function $\psi:X\rightarrow\mathbb{R}$ satisfying 
    $$|\varphi(x)-\psi(x)|<\varepsilon\quad \textnormal{for all } x\in K \textnormal{ and } \psi|_A=0.$$
    
    We may assume that $X$ is a Zariski closed subvariety of $\mathbb{R}^n$ for some $n$. Then $A=q^{-1}(0)$ for some polynomial function $q:\mathbb{R}^n\rightarrow \mathbb{R}$. Note that the map 
    $$e:X\rightarrow \mathbb{R}^{n+1},\quad x=(x_1,\ldots, x_n)\mapsto (q(x),x_1,\ldots,x_n)$$
    induces a biregular isomorphism from $X$ onto the Zariski closed subvariety $e(X)$ of $\mathbb{R}^{n+1}$ with 
    $$e(A)=e(X)\cap H,$$
    where $H=\{ (x_0,\ldots,x_n)\in\mathbb{R}^{n+1}:\, x_0=0 \}$.
    
    According to Tietze's extension theorem, there exists a continuous function $\alpha:\mathbb{R}^{n+1}\rightarrow\mathbb{R}$ such that 
    $$\alpha\circ e|_K=\varphi|_K\,\textnormal{ and }\,\alpha|_H=0.$$
    In turn, by Whitney's approximation theorem \cite[(15.8.1)]{22} there exists a smooth function $\beta:\mathbb{R}^{n+1}\rightarrow\mathbb{R}$ satisfying
    $$|\alpha(y)-\beta(y)|<\frac{\varepsilon}{2}\quad \textnormal{for all } y\in e(K) \textnormal{ and } \beta|_H=0.$$
    The function $\beta$, being smooth, can be written as $b=\pi\gamma$, where $\pi:\mathbb{R}^{n+1}\rightarrow\mathbb{R}$ is the projection $\pi(y_0,\ldots,y_n)=y_0$ and $\gamma:\mathbb{R}^{n+1}\rightarrow\mathbb{R}$ is a smooth function. By the Stone-Weierstrass approximation theorem, there exists a polynomial function $p:\mathbb{R}^{n+1}\rightarrow\mathbb{R}$ such that 
    $$|\gamma(y)-p(y)|<\frac{\varepsilon}{2M}\quad \textnormal{ for all } y\in e(K),$$
    where $M:=\textnormal{sup}\{ |\pi(y)|:\, y\in e(K) \}$. Therefore
    $$|\beta(y)-\pi(y)p(y)|=|\pi(y)\gamma(y)-\pi(y)p(y)|<\frac{\varepsilon}{2}\quad \textnormal{ for all } y\in e(K).$$
    Finally,
    $$|\alpha(y)-\pi(y)p(y)|\leq |\alpha(y)-\beta(y)|+|\beta(y)-\pi(y)p(y)|<\frac{\varepsilon}{2}+\frac{\varepsilon}{2}=\varepsilon$$
    for all $y\in e(K)$, so the regular function $\psi:=(\pi\circ e)(p\circ e)$ on $X$ has the required properties.
\end{proof}

It readily follows that Proposition \ref{prop:4.4} is equivalent to Theorem \ref{thm:1.14} with $Y=\mathbb{R}$ (see the comment after the proof of Proposition \ref{prop:4.1}).

\section{Real algebraic varieties with the AIP}
\label{section:5}

In this section, we prove Theorems \ref{thm:1.6} and \ref{thm:1.14}. Our first result below is stronger than the implication (e) $\Rightarrow$ (a) in Theorem \ref{thm:1.6}.

\begin{theorem}
\label{thm:5.1}
    Let $Y$ be a malleable nonsingular real algebraic variety. Let $X$ be a nonsingular real algebraic variety, $f:U\rightarrow Y$ a smooth map defined on an open subset $U$ of $X$, and $A$ a Zariski closed subvariety of $X$ which is faithful at every point in $U\cap A$. Assume that there exists a regular map $f_0:X\rightarrow Y$ with $f_0|_{U\cap A}=f|_{U\cap A}$ and with $f_0|_U$ homotopic to $f$ relative to $U\cap A$. Then $f$ can be approximated in the $\mathcal{C}^{\infty}$ topology by regular maps from $X$ to $Y$ whose restrictions to $A$ are equal to $f_0|_A$.
\end{theorem}
\begin{proof}
    In the proof, all statements about approximation refer to approximation in the $\mathcal{C}^{\infty}$ topology. It is sufficient to prove that for any given open set $U_0\in\Omega(U)$, the restriction $f|_{U_0}$ can by approximated by regular maps from $X$ to $Y$ whose restrictions to $A$ are equal to $f_0|_A$.

    By assumption, there exists a continuous map $F:U\times [0,1]\rightarrow Y$ such that $F_0=f_0|_U$, $F_1=f$ and $F_t|_{U\cap A}=f_0|_{U\cap A}$ for all $t\in [0,1]$. We may assume that $F$ is a smooth map \cite[Proposition 10.22]{44}. Consider the subset 
    \begin{align*}
        I:=\{ t\in [0,1]:\, F_t|_{U_0} \textnormal{ can be approximated by regular maps from } X \textnormal{ to } Y&\\ 
        \textnormal{ whose restrictions to } A  \textnormal{ are equal to } f_0|_A \}&
    \end{align*}
    of $[0,1]$. By definition, $I$ is closed and $0\in I$. We now prove that $I$ is open in $[0,1]$. Let $t_0\in I$. Let dist be a metric on $Y$. According to \cite[Lemma 2.2]{17}, $Y$ admits a dominating spray $(E,p,s)$, where $p:E=Y\times \mathbb{R}^n\rightarrow Y$  is the product vector bundle. Choose a constant $\varepsilon>0$ as in Lemma \ref{lem:2.5} (with $N=Y$, $M=U$, $B=U_0\cap A$, $K=\bar{U}_0$, $f=F_{t_0}|_K$). Now choose a neighborhood $I_0$ of $t_0$ in $[0,1]$ such that 
    $$\textnormal{dist}(F_t(x),F_{t_0}(x))<\varepsilon\quad \textnormal{for all } t\in I_0 \textnormal{ and all } x\in K.$$
    Thus, for each $t\in I_0$ there exists a continuous map $\xi:K\rightarrow E=Y\times \mathbb{R}^n$ satisfying 
    \begin{enumerate}
        \item[(i)] $p\circ\xi = F_{t_0}|_K$,
        \item[(ii)] $s\circ\xi = F_t|_K$,
        \item[(iii)] $\xi|_{U_0}$ is smooth,
        \item[(iv)] $\xi(x)=0_{F_{t_0}(x)}$ for all $x\in U_0\cap A$.
    \end{enumerate}
    It follows from (i) that the map $\xi$ is of the form $\xi=(F_{t_0}|_K,\eta)$ for some continuous map $\eta:K\rightarrow\mathbb{R}^n$. The restriction $\eta|_{U_0}$ is a smooth map by (iii), and $\eta|_{U_0\cap A}=0$ by (iv). Moreover, in view of (ii) we get 
    $$F_t(x)=s(F_{t_0}(x),\eta(x))\quad \textnormal{ for all } x\in U_0.$$
    According to the definition of $I$, the map $F_{t_0}|_{U_0}$ can be approximated by regular maps from $X$ to $Y$ whose restricitons to $A$ are equal to $f_0|_A$. By Proposition \ref{prop:4.1}, the map $\eta$ can be approximated by regular maps from $X$ to $\mathbb{R}^n$ whose restrictions to $A$ are identically $0$. Thus, since $s$ is a regular map, the map $F_t|_{U_0}$ can be approximated by regular maps from $X$ to $Y$ whose restrictions to $A$ are equal to $f_0|_A$, so $t\in I$. Therefore the inclusion $I_0\subset I$ holds, which means that $I$ is an open subset of $[0,1]$. In summary, $I=[0,1]$ due to the connectedness of the interval $[0,1]$. The proof is complete because $f|_{U_0}=F_1|_{U_0}$.
\end{proof}

The following result will be useful in the proof of Theorem \ref{thm:1.6}.

\begin{theorem}
\label{thm:5.2}
    Let $Y$ be a nonsingular real algebraic variety. Let $(E,p,\sigma)$ be a dominating smooth spray for $Y$, where $p:E\rightarrow Y$ is an algebraic vector bundle over $Y$. Assume that the smooth map $\sigma :E\rightarrow Y$ can be approximated in the $\mathcal{C}^1$ topology by regular maps from $E$ to $Y$ whose restrictions to the zero section $Z(E)$ are equal to $p|_{Z(E)}$. Then $Y$ is a malleable variety.
\end{theorem}
\begin{proof}
    By \cite[Theorem 12.1.7]{5}, there exists an algebraic vector bundle $E'$ over $Y$ such that the direct sum $E\oplus E'$ is algebraically isomorphic to the product vector bundle $Y\times \mathbb{R}^n$. We identify $Y\times \mathbb{R}^n$ with $E\oplus E'$ and denote by $\pi:Y\times \mathbb{R}^n\rightarrow E$ the projection morphism. Note that $(Y\times \mathbb{R}^n,p\circ\pi,\sigma\circ\pi)$ is a dominating smooth spray for $Y$, and the smooth map $\sigma\circ\pi:Y\times \mathbb{R}^n\rightarrow Y$ can be approximated in the $\mathcal{C}^1$ topology by regular maps from $Y\times \mathbb{R}^n$ to $Y$ whose restrictions to $Z(Y\times \mathbb{R}^n)=Y\times \{0\}$ are equal to $p\circ \pi|_{Y\times \{0\}}$. This reduces the proof to the case where $p:E=Y\times \mathbb{R}^n\rightarrow Y$ is the product vector bundle.

    Using the assumption of approximability of $\sigma$, for each point $b\in Y$ we can find a regular map $s_b:Y\times \mathbb{R}^n\rightarrow Y$ and a Zariski open neighborhood $V_b\subset Y$ of $b$ such that $s_b|_{Y\times \{0\}}=p|_{Y\times \{0\}}$ and for every point $y\in V_b$ the regular map 
    $$\mathbb{R}^n\rightarrow Y,\quad v\mapsto s_b(y,v)$$
    is a submersion at $0\in\mathbb{R}^n$.
    Since $Y$ is a compact space in the Zariski topology, we get $Y=V_{b_1}\cup\ldots\cup V_{b_r}$ for some points $b_1$, \ldots, $b_r$. Set $V_i:=V_{b_i}$ and $s_i:=s_{b_i}$. For $i\in \{ 1,\ldots, r \}$, we define recursively a regular map $s^{(i)}:Y\times (\mathbb{R}^n)^i\rightarrow Y$ by
    $$s^{(1)}=s_1,$$
    $$s^{(i)}(y,v_1,\ldots, v_i)=s_i(s^{(i-1)}(y,v_1,\ldots,v_{i-1}),v_i)\quad \textnormal{for } i\geq 2,$$
    where $v_1$, \ldots, $v_i$ belong to $\mathbb{R}^n$. Let $p^{(r)}:E^{(r)}=Y\times (\mathbb{R}^n)^r\rightarrow Y$ be the product vector bundle. Then $(E^{(r)},p^{(r)},s^{(r)})$ is a dominating spray for $Y$, so $Y$ is a malleable variety.
\end{proof}

We are now ready to perform the main tasks of this section.

\begin{proof}[Proof of Theorem \ref{thm:1.6}]
    Obviously, (a) implies (b). In light of Observation \ref{obs:3.4}, (b) implies both (c) and (d). By Theorem \ref{thm:5.1}, (e) implies (a).

    According to Lemma \ref{lem:2.8} and Lemma \ref{lem:2.3} (with $E^{\circ}=E$), there exists a dominating smooth spray for $Y$ of the form $(TY,p,\sigma)$, where $p:TY\rightarrow Y$ is the tangent bundle of $Y$. Therefore (c) implies (e) by Theorem \ref{thm:5.2}.

    Suppose (d) holds. Then, by Lemma \ref{lem:2.8} there exists a dominating smooth spray $(E,p,\sigma)$ for $Y$, where $p:E=Y\times \mathbb{R}^n\rightarrow Y$ is the product vector bundle and the smooth map $\sigma:E\rightarrow Y$ can be approximated in the $\mathcal{C}^1$ topology by regular maps from $E$ to $Y$ whose restrictions to $Z(E)=Y\times \{0\}$ are equal to $p|_{Y\times \{0\}}$. Hence (e) holds by Theorem \ref{thm:5.2}.
\end{proof}
\begin{proof}[Proof of Theorem \ref{thm:1.14}]
    In the proof, all statements about approximation refer to approximation in the compact-open topology. Our goal is to show that for any given compact subset $K$ of $D$, the restriction $f|_K$ can be approximated by regular maps from $X$ to $Y$ whose restrictions to $A$ are equal to $f_0|_A$. We may assume that $X$ is a Zariski closed subvariety of $\mathbb{R}^m$ for some $m$.

    By assumption, there exists a continuous map $F:D\times [0,1]\rightarrow Y$ such that $F_0=f_0|_D$, $F_1=f$ and $F_t|_{D\cap A}=f_0|_{D\cap A}$ for all $t\in [0,1]$. Consider the subset 
    \begin{align*}
        I:=\{ t\in [0,1]:\, F_t|_K \textnormal{ can be approximated by regular maps from } X \textnormal{ to } Y&\\ 
        \textnormal{ whose restrictions to } A \textnormal{ are equal to } f_0|_A  \}&
    \end{align*}
    of $[0,1]$. By definition, $I$ is closed and $0\in I$. We now prove that $I$ is open in $[0,1]$. Let $t_0\in I$. Let dist be a metric on $Y$. According to \cite[Lemma 2.2]{17} $Y$ admits a dominating spray $(E,p,s)$, where $p:E=Y\times \mathbb{R}^n\rightarrow Y$ is the product vector bundle over $Y$. Choose a constant $\varepsilon>0$ as in Lemma \ref{lem:2.5} (with $N=Y$, $M=\mathbb{R}^m$, $M_0=\varnothing$, $B=K\cap A$, $f=F_{t_0}|_K$). Now choose a neighborhood $I_0$ of $t_0$ in $[0,1]$ such that 
    $$\textnormal{dist}(F_t(x),F_{t_0}(x))<\varepsilon\quad \textnormal{for all } t\in I_0 \textnormal{ and all } x\in K.$$
    Thus, for each $t\in I_0$ there exists a continuous map $\xi: K\rightarrow E=Y\times \mathbb{R}^n$ satisfying 
    \begin{enumerate}
        \item[(i)] $p\circ\xi=F_{t_0}|_K$,
        \item[(ii)] $s\circ\xi=F_t|_K$,
        \item[(iii)] $\xi(x)=0_{F_{t_0}(x)}$ for all $x\in K$.
    \end{enumerate}
     It follows from (i) that the map $\xi$ is of the form $\xi=(F_{t_0},\eta)$ for some continuous map $\eta:K\rightarrow\mathbb{R}^n$. By (iii), the map $\eta$ satisfies $\eta|_{K\cap A}=0$. Moreover, in view of (ii), we get
     $$F_t(x)=s(F_{t_0}(x),\eta(x))\quad \textnormal{for all } x\in K.$$
     According to the definition of $I$, the map $F_{t_0}|_K$ can be approximated by regular maps from $X$ to $Y$ whose restrictions to $A$ are equal to $f_0|_A$. By Proposition \ref{prop:4.4}, the map $\eta$ can be approximated by regular maps from $X$ to $\mathbb{R}^n$ whose restrictions to $A$ are identically $0$. Thus, since $s$ is a regular map, the map $F_t|_K$ can be approximated by regular maps from $X$ to $Y$ whose restrictions to $A$ are equal to $f_0|_A$, so $t\in I$. Therefore the inclusion $I_0\subset I$ holds, which means that $I$ is an open subset of $[0,1]$. In summary, $I=[0,1]$ due to the connectedness of the interval $[0,1]$. The proof is complete because $f|_K=F_1|_K$.
\end{proof}

\section{Extensions of local sprays}
\label{section:6}

In this section, we introduce and study local sprays for nonsingular real algebraic varieties. The resulting Theorem \ref{thm:6.3} directly implies Theorem \ref{thm:1.9}.

\begin{definition}
\label{def:6.1}
    Let $Y$ be a nonsingular real algebraic variety.
    \begin{itemize}
        \item[(1)] (Local spray).  A \emph{local spray} for $Y$ is a regular map $\sigma:U\times\mathbb{R}^n\rightarrow Y$, where $U$ is a Zariski open subset of $Y$ and $n$ is a nonnegative integer, such that $\sigma(y,0)=0$ for all $y\in U$.
        \item[(2)] (Dominating local spray).  A local spray $\sigma:U\times\mathbb{R}^n\rightarrow Y$ for $Y$ is said to be \emph{dominating} if for every point $y\in U$ the regular map 
        $$\mathbb{R}^n\rightarrow Y,\quad v\mapsto\sigma(y,v)$$
        is a submersion at $0\in\mathbb{R}^n$.
        \item[(3)] (Locally malleable variety).  The variety $Y$ is called \emph{locally malleable} if for every point $y_0\in Y$ there exists a dominating local spray $\sigma:U\times \mathbb{R}^n\rightarrow U$ for $Y$ with $y_0\in U$.
    \end{itemize}
\end{definition}

Definition \ref{def:6.1}(2) and Definition \ref{def:1.2}(6) (with $Y=U$) are interrelated as follows.

\begin{proposition}
\label{prop:6.2}
    Let $Y$ be a nonsingular real algebraic variety, and $U$ a Zariski open subset of $Y$. Then the following conditions are equivalent:
    \begin{enumerate}
        \item[(a)] There exists a dominating local spray $\sigma:U\times\mathbb{R}^n\rightarrow Y$ for $Y$.
        \item[(b)] The variety $U$ is malleable.
    \end{enumerate}
\end{proposition}
\begin{proof}
    Suppose (a) holds. Then $\sigma^{-1}(U)$ is a Zariski open neighborhood of $U\times \{0\}$ in $U\times\mathbb{R}^n$, so by \cite[Lemma 2.10]{41} there exists a regular function $\varepsilon:U\rightarrow \mathbb{R}$ such that 
    $$\varepsilon(y)>0\textnormal{ and } (y,\varepsilon(y)\frac{v}{1+||v||^2})\in\sigma^{-1}(U)\textnormal{ for all } (y,v)\in U\times\mathbb{R}^n,$$
    where $||\cdot||$ is the Euclidean norm on $\mathbb{R}^n$. We obtain a dominating spray $(E,p,s)$ for $U$, where $p:E=U\times\mathbb{R}^n\rightarrow U$ is the product vector bundle over $U$ and $s:E\rightarrow U$ is defined by 
    $$s(y,v)=\sigma(y,\varepsilon(y)\frac{v}{1+||v||^2})\textnormal{ for all } (y,v)\in U\times\mathbb{R}^n.$$
    Therefore $U$ is a malleable variety.

    On the other hand, if (b) holds, then $U$ admits a dominating spray $(E,p,s)$ where $p:E=U\times\mathbb{R}^n\rightarrow U$ is the product vector bundle over $U$ \cite[Lemma 2.2]{17}. Obviously, $s:E\rightarrow U\subset Y$ can be regarded as a dominating local spray for $Y$, which gives (a).
\end{proof}

Here is the main result of this section.

\begin{theorem}
\label{thm:6.3}
    Let $Y$ be a nonsingular real algebraic variety that is locally malleable. Then $Y$ is malleable.
\end{theorem}
\begin{proof}
    Since $Y$ is a compact space in the Zariski topology, there exists a finite collection $\{ \sigma_i:U_i\times\mathbb{R}^{n_i}\rightarrow Y\,:\, i=1,\ldots, r \}$ of dominating local sprays for $Y$, where the Zariski open sets $U_i$ cover $Y$. By Lemma \ref{lem:6.4} below, for each $i\in\{1,\ldots, r\}$ there exists a regular function $\varphi_i:Y\rightarrow\mathbb{R}$ such that $\varphi_i^{-1}(0)=Y\setminus U_i$ and the map $\tau_i:Y\times \mathbb{R}^{n_i}\rightarrow Y$ defined by
    $$\tau_i(y,v_i):=
    \begin{cases}
        \sigma_i(y,\varphi_i(y)v_i) \textnormal{ for } (y,v_i)\in U_i\times\mathbb{R}^{n_i}\\
        y\qquad\qquad\quad\, \textnormal{ for } (y,v_i)\in(Y\setminus U_i)\times\mathbb{R}^{n_i}
    \end{cases}$$
    is regular. By construction, for every point $y\in U_i$ the regular map
    $$\mathbb{R}^{n_i}\rightarrow Y,\quad v_i\mapsto \tau_i(y,v_i)$$
    is a submersion at $0\in\mathbb{R}^{n_i}$. For $i\in \{1,\ldots,r\}$, we define recursively a regular map
    $$s_i:Y\times\mathbb{R}^{n_1}\times \ldots \times \mathbb{R}^{n_i} \rightarrow Y$$ 
    by
    $$s_1=\tau_1,$$
    $$s_i(y,v_1,\ldots,v_{i-1},v_i)=\tau_i(s_{i-1}(y,v_1,\ldots,v_{i-1}),v_i)\, \textnormal{ for } i\geq 2.$$
    We obtain a dominating spray $(E,p,s)$ for $Y$, where
    $$p:E=Y\times\mathbb{R}^{n_1}\times\ldots\times \mathbb{R}^{n_r}\rightarrow Y$$
    is the product vector bundle and $s=s_r$. Thus, $Y$ is a malleable variety.
\end{proof}

In the proof of Theorem \ref{thm:6.3}, the key role was played by the following extension lemma for local sprays.

\begin{lemma}
\label{lem:6.4}
    Let $Y$ be a nonsingular real algebraic variety, $U$ a Zariski open subset of $Y$, and $\sigma:U\times \mathbb{R}^n\rightarrow Y$ a regular map satisfying $\sigma(y,0)=y$ for all $y\in U$. Then there exists a regular function $\varphi:Y\rightarrow\mathbb{R}$ with $\varphi^{-1}(0)=Y\setminus U$ and such that the map $\tau:Y\times \mathbb{R}^n\rightarrow Y$ defined by
    $$\tau(y,v):=\begin{cases}
        \sigma(y,\varphi(y)v)\textnormal{ for } (y,v)\in U\times \mathbb{R}^n\\
        y\qquad\qquad\,\,\, \textnormal{for } (y,v)\in(Y\setminus U)\times \mathbb{R}^n
    \end{cases}$$
    is regular.
\end{lemma}
\begin{proof}
    Note that in the ring of regular functions on $U\times\mathbb{R}^n$, the ideal of all functions vanishing on $U\times \{0\}\subset U\times \mathbb{R}^n$ is generated by the projections 
    $$U\times \mathbb{R}^n\rightarrow \mathbb{R},\quad (y,v)\mapsto v_j$$
    where $v=(v_1,\ldots, v_n)$ and $j=1,\ldots,n$. The analogous statement is also true if $U$ is replaced by $Y$.

    We may assume that $Y$ is an algebraic subset of $\mathbb{R}^m$. Then 
    $$\sigma(y,v)=(\sigma_1(y,v),\ldots,\sigma_m(y,v))\textnormal{ for all } (y,v)\in U\times\mathbb{R}^n,$$
    where each $\sigma_i$ is a regular function on $U\times \mathbb{R}^n$ for $i=1,\ldots, m$. Moreover,
    $$\sigma_i(y,0)=y_i\textnormal{ for all } y=(y_1,\ldots,y_m)\in U.$$
    By the observation at the beginning of the proof, $\sigma_i(y,v)$ can be written as 
    $$\sigma_i(y,v)=y_i+\sum_{j=1}^n \alpha_{ij}(y,v)v_j$$
    for some regular functions $\alpha_{ij}$ on $U\times\mathbb{R}^n$. In turn, each $\alpha_{ij}$ can be written as the quotient of two regular functions on $Y\times \mathbb{R}^n$, where the denominator does not vanish at any point of $U\times\mathbb{R}^n$ (see \cite[p.62]{5} or \cite[p.14]{49}). Taking the common denominator $\gamma$ of such quotients for all $i$ and $j$, we get $\gamma^{-1}(0)\subset (Y\setminus U)\times\mathbb{R}^n$ and 
    $$\alpha_{ij}(y,v)=\frac{\beta_{ij}(y,v)}{\gamma(y,v)}\textnormal{ for all } (y,v)\in U\times\mathbb{R}^n$$
    for some regular functions $\beta_{ij}$ on $Y\times \mathbb{R}^n$. Again, using the observation at the beginning of the proof, $\gamma(y,v)$ can be written as 
    $$\gamma(y,v)=\gamma(y,0)+\sum_{j=1}^n\delta_j(y,v)v_j\textnormal{ for all } (y,v)\in Y\times\mathbb{R}^n$$
    for some regular functions $\delta_j$ on $Y\times \mathbb{R}^n$.

    Let $\psi$ be any regular function on $Y$ with $\psi^{-1}(0)=Y\setminus U$. Define a regular function $\varphi:Y\rightarrow\mathbb{R}$ by
    $$\varphi(y)=\psi(y)\gamma(y,0) \textnormal{ for all } y\in Y.$$
    Then $\varphi^{-1}(0)=Y\setminus U$. Moreover,
    \begin{align*}
        \sigma_i(y,\varphi(y)v)&=
        y_i+\frac{\sum_{j=1}^n \beta_{ij}(y,\varphi(y)v)\varphi(y)v_j}{\gamma(y,\varphi(y)v)}\\
        &= y_i+ \frac{\sum_{j=1}^n \beta_{ij}(y,\varphi(y)v)\psi(y)\gamma(y,0)v_j}{\gamma(y,0)+\sum_{j=1}^n \delta_j(y,\varphi(y)v)\psi(y)\gamma(y,0)v_j}\\
        &= y_i+ \frac{\psi(y)\sum_{j=1}^n \beta_{ij}(y,\varphi(y)v)v_j}{1+\psi(y)\sum_{j=1}^n \delta_j(y,\varphi(y)v)v_j}.
    \end{align*}
    Since $\gamma(y,v)\neq 0$ for all $(y,v)\in U\times \mathbb{R}^n$, the denominator 
    $$d_i(y,v):=1+\psi(y)\sum_{j=1}^n \delta_j(y,\varphi(y)v)v_j$$
    has the same property. On the other hand, $d_i(y,v)=1$ for all $(y,v)\in(Y\setminus U)\times\mathbb{R}^n$. It follows that $\tau:Y\times \mathbb{R}^n\rightarrow Y$ is a regular map, as required.
\end{proof}

The first attempt to prove Theorem \ref{thm:6.3} was made in \cite{41}, but only led to a weaker version \cite[Proposition 2.8]{41}, because Lemma \ref{lem:6.4} was not available in its current strong form.

\paragraph{Acknowledgements.}We are grateful to O. Benoist and O. Wittenberg for clarifying some issues regarding the rationality of homogeneous spaces, and to J. Bochnak for questions and comments. Both authors were partially supported by the National Science Center (Poland) under grant number 2022/47/B/ST1/00211.

\vspace{20pt}

\noindent
Institute of Mathematics\\
Faculty of Mathematics and Computer Science\\
Jagiellonian University\\
Łojasiewicza 6\\
30-348 Kraków\\
Poland\\
E-mail: Juliusz.Banecki@student.uj.edu.pl\\
E-mail: Wojciech.Kucharz@im.uj.edu.pl

\end{document}